\documentclass{amsart}

\usepackage{amsmath}
\usepackage{amscd}
\usepackage{amssymb} 

\newcommand{\cal}{\mathcal}
\newcommand{\bk}{{\bf k}}

\newcommand{\bC}{{\Bbb C}}

\newcommand{\bQ}{{\Bbb Q}}
\newcommand{\bR}{{\Bbb R}}

\newcommand{\cA}{{\cal A}}

\newcommand{\cH}{{\cal H}}
\newcommand{\cL}{{\cal L}}

\newcommand{\fk}{{\frak k}}

\newcommand{\ft}{{\frak t}}

\DeclareMathOperator{\Img}{Im}

\DeclareMathOperator{\Ker}{Ker}

\newtheorem{theorem}{Theorem}[section]
\newtheorem{theorem/definition}{Theorem/Definition}[section]

\newtheorem{proposition}{Proposition}[section]
\newtheorem{lemma}{Lemma}[section]

\newtheorem{corollary}{Corollary}[section]

\theoremstyle{remark}
\newtheorem{remark}{Remark}[section]

\theoremstyle{definition}

\begin{document}
\title
{Formal Frobenius manifold structure on equivariant cohomology}
\author{Huai-Dong Cao, Jian Zhou}
\address{Department of Mathematics\\
Texas A\&M University\\
College Station, TX 77843}
\email{cao@math.tamu.edu, zhou@math.tamu.edu}
\begin{abstract}
For a closed K\"{a}hler manifold with a Hamiltonian action 
of a connected compact Lie group by holomorphic isometries,
we construct a formal Frobenius manifold structure on the equivariant cohomology
by exploiting a natural DGBV algebra structure on the Cartan model.
\end{abstract}
\maketitle

The notion of Frobenius manifolds was introduced 
by Dubrovin \cite{Dub1, Dub2}.
It gives a coordinate free formulation of solutions to the WDVV equations.
As surveyed in Manin \cite{Man2},
there are three major methods to construct solutions to WDVV equations.
The first method involves
the theory of quantum cohomology via Gromov-Witten invariants
(or topological sigma model in physics literature),
see e.g. Ruan-Tian \cite{Rua-Tia} and Kontsevich-Manin \cite{Kon-Man}.
The second method is Saito's theory of singularities 
(or Landau-Ginzburg model in physics literature).
The third method exploits the so-called DGBV algebras,
named after Gerstenhaber, Batalin and Vilkovisky.
This last method first appeared in Barannikov-Kontsevich \cite{Bar-Kon} in the context of extended 
moduli spaces of Calabi-Yau manifolds,  based on the Kodaira-Spencer theory of gravity 
of Bershadsky-Cecotti-Ooguri-Vafa \cite{Ber-Cec-Oog-Vaf}
which extends earlier works of Tian \cite{Tia} and Todorov \cite{Tod}.
A detailed account of the construction for general DGBV algebras can be found in 
Manin \cite{Man2}.

GBV algebras have appeared in many places in Mathematics and
Mathematical Physics,
e.g. algebraic deformation theory and Hochschild cohomology 
(Gerstenhaber \cite{Ger}),
string theory (Lian-Zuckerman \cite{Lia-Zuc}),
gauge theory (Batalin-Vilkovisky \cite{Bat-Vil}), etc.
However,  examples of DGBV algebras in differential geometry 
were relatively rare.
Earlier examples include Tian's formula \cite{Tia} in 
deformation theory of Calabi-Yau manifolds and
Koszul's operator in Poisson geometry \cite{Kos}.
But the recognizations of DGBV algebra structures in these theories seem to come later
in e.g. Ran \cite{Ran} and Xu \cite{Xu} respectively.
In a series of papers \cite{Cao-Zho1, Cao-Zho2, Cao-Zho3},
the authors constructed many DGBV algebras from K\"{a}hler 
and hyperk\"{a}hler manifolds,
and showed that they satisfy the conditions to carry out 
the construction of formal Frobenius manifold structures on the cohomology.
Also, it was shown that different DGBV algebra structures 
can yield the same solution to the WDVV equations.
In particular, we get formal Frobenius manifold structures on the de Rham and
Dolbeault cohomology of a closed K\"{a}hler manifold.
In this paper, we carry over the same ideas to equivariant cohomology,
in the case 
of closed K\"{a}hler manifolds with Hamiltonian actions of a Lie
group by holomorphic isometries.

The main result in this paper is related to the equivariant quantum cohomology.
There are three models to define equivariant cohomology:
the Borel model, the Cartan model and the Weil model. 
In Givental-Kim \cite{Giv-Kim},
a version of quantum cohomology based on Borel model was suggested and
the rigorous formulation appeared in Lu \cite{Lu}. 
Some discussions of WDVV equations and Frobenius manifold structure in 
equivariant quantum cohomology can be found in Givental \cite{Giv}.
Our construction of a formal Frobenius manifold structure on equivariant cohomology
uses the Cartan model,
which enables us to manipulate everything by differential forms. 
Differential geometers are familiar with the idea that 
choosing nice representatives of  cohomology classes by differential forms
may lead to  more information.
This idea proves useful again in our theory:
as shown in \S 6, we can work over the ring
given by the equivariant cohomology of a point,
while Givental \cite{Giv} has to use its fractional field.

The rest of the paper is arranged as follows. 
We review some definitions and the construction of 
formal Frobenius manifolds from DGBV algebras in \S 1.
The Cartan model for equivariant cohomology is reviewed in \S 2.
In \S 3, we construct a DGBV algebra structure on the Cartan model
when the group preserves a Poisson structure.
We begin in \S 4 the discussion of Hamiltonian actions.
The main results appear in the more technical sections \S 5 and \S 6,
where we restrict our attention to the K\"{a}hler case. 

{\bf Acknowledgements}. 
{\em  The work in this paper is carried out
while the second author is visiting Texas A$\&$M University.
He thanks the Mathematics Department 
and the Geometry-Analysis-Topology group for the hospitality
and financial support.}

\section{DGBV algebras and formal Frobenius manifolds}
\label{sec:DGBV}

\subsection{Frobenius algebra}

Let $\bk$ be a commutative $\bQ$-algebra,
$H$ a free $\bk$-module of finite rank, 
endowed with a commutative associative multiplication
$$\wedge: \;\; H \otimes_{\bk} H \to H.$$
We call $(H, \wedge)$ a {\em Frobenius algebra} 
if there is a symmetric nondegenerate bilinear form
$(\cdot, \cdot): H \otimes H \to \bk$ 
such that
\begin{eqnarray} \label{eqn:Frobenius}
(X \wedge Y, Z) = (X, Y \wedge Z)
\end{eqnarray}
for any $X, Y, Z \in H$.
Such a bilinear form $(\cdot, \cdot)$ is called an {\em invariant inner product}
on $H$.

Take a basis $\{ e_a \}$ of $H$.
Let $\eta_{ab} = (e_a, e_b)$ and $(\eta^{ab})$ 
be the inverse matrix of $(\eta_{ab})$.
Also let $\phi_{ab}^c$ be the structure constants defined by
$$e_a \wedge e_b = \phi_{ab}^c e_c.$$
It is clear that the constants $\eta_{ab}$ and $\phi_{ab}^c$'s completely determine the structure
of the Frobenius algebra.
When $(H, \wedge)$ has an identity $1$, $\eta_{ab}$ and $\phi_{ab}^c$'s
 can be encoded in a symmetric $3$-tensor $\phi \in S^3H^*$ as follows.
Assume that $e_0 = 1$.
Set $\phi_{abc} = \phi_{ab}^p\eta_{pc}$.
Then 
$$\phi_{abc} = (e_a \wedge e_b, e_c).$$
From (\ref{eqn:Frobenius}),
one sees that $\phi$ is symmetric in the three indices.
One can recover the inner product and the multiplication from $\phi$ since
\begin{eqnarray*}
\eta_{ab} = \phi_{0ab}, & \phi_{ab}^c = \phi_{abp}\eta^{pc}.
\end{eqnarray*}
The associativity of the multiplication is equivalent to the following
system of equations
\begin{eqnarray} \label{eqn:WDVV1}
\phi_{abp}\eta^{pq}\phi_{qcd} = \phi_{bcp}\eta^{pq}\phi_{aqd}.
\end{eqnarray}

\subsection{WDVV equations and Frobenius manifolds}
Let $(H, \wedge, (\cdot, \cdot))$ be a finite dimensional 
Frobenius algebra with $1$ over $\bk$.
Let $\{e_a \}$ be a basis of $H$ as above.
Denote by $\{x^a\}$ the linear coordinates in the basis $\{e_a\}$.
Consider a self-parameterizing family $(H, \{\cdot_x, x \in U\})$, where $U$
is an open subset of $H$,
such that $1$ is the identity for each $\wedge_x$ and
$$(X \wedge_x Y, Z) = (X, Y \wedge_x Z),$$
for all $X, Y, Z \in H$ and $x \in U$.
Then we get a family of $3$-tensors $\phi_{abc}(x)$.
String theory (see e.g. Dijkgraaf-Verlinde-Verlinde \cite{Dij-Ver-Ver}) 
suggests that one should require
$$\frac{\partial}{\partial x^d} \phi_{abc} 
= \frac{\partial}{\partial x^c} \phi_{abd}.$$
Under this condition,
if $U$ is contractible,
one can find a function $\Phi: U \to \bk$,
such that
$$\phi_{abc}(x) = \frac{\partial^3\Phi}{\partial x^a \partial x^b\partial x^c}.$$
By (\ref{eqn:WDVV1}), the associativity of $\wedge_x$ is then equivalent to that
$\Phi$ satisfies the following Witten-Dijkgraaf-E. Verlinde-H. Verlinde
(WDVV) equations:
\begin{eqnarray} \label{eqn:WDVV}
\frac{\partial^3 \Phi}{\partial x^a \partial x^b \partial x^p} \eta^{pq}
\frac{\partial^3 \Phi}{\partial x^q \partial x^c \partial x^d} =
\frac{\partial^3 \Phi}{\partial x^b \partial x^c \partial x^p} \eta^{pq}
\frac{\partial^3 \Phi}{\partial x^a \partial x^q \partial x^d}.
\end{eqnarray}
Such a function $\Phi$ is called a {\em potential function}.
Dubrovin \cite{Dub1, Dub2} introduced the notion of a Frobenius manifold
to give a global formulation.
For our purpose in this paper, a Frobenius manifold structure on 
$(H, \wedge, (\cdot, \cdot))$ will mean a 
solution $\Phi$ to the WDVV equations with
\begin{eqnarray} \label{eqn:metric}
\eta_{ab} 
= \frac{\partial^3\Phi}{\partial x^0\partial x^a \partial x^b}.
\end{eqnarray}
By definition a structure of {\em formal Frobenius manifold} on
$(H, \wedge, (\cdot, \cdot))$ is a formal power series $\Phi$ 
which satisfies the WDVV equations.
We refer to $(H, \wedge, (\cdot, \cdot))$
as the initial data for the WDVV equations.
If $\Phi$ also satisfies (\ref{eqn:metric}),
it is called a structure of formal Frobenius manifold with identity.
The above discussion can be carried out for graded algebras as well.
See Manin \cite{Man1}.

\subsection{DGBV algebras and Frobenius manifolds}

Let $(\cA, \wedge)$ be a graded 
commutative associative algebra over ${\bk}$. 
For any linear operator $\Delta$ of odd degree,
define
$$[a \bullet b]_{\Delta} = (-1)^{|a|} (\Delta (a \wedge b) 
- (\Delta a) \wedge b - (-1)^{|a|} a \wedge \Delta b),$$
for homogeneous elements $a, b \in \cA$.
If $\Delta^2 = 0$ and
\begin{eqnarray*}
[a \bullet (b \wedge c)]_{\Delta} 
	= [a \bullet b]_{\Delta} \wedge c
	+ (-1)^{(|a|+1)|b|} b \wedge [a \bullet c]_{\Delta},
\end{eqnarray*}
for all homogeneous $a, b, c \in \cA$,
then $(\cA, \wedge, \Delta, [\cdot\bullet\cdot]_{\Delta})$ is 
a {\em Gerstenhaber-Batalin-Vilkovisky (GBV) algebra}.
(Notice that if one takes $a=b=c=1$,
then one can deduce $\Delta 1 = 0$.)
A  {\em DGBV (differential Gerstenhaber-Batalin-Vilkovisky) algebra}
is a GBV algebra with a $\bk$-linear derivation $\delta$ of odd degree
with respect to $\wedge$, such that
$$ \delta^2 = \delta\Delta + \Delta \delta = 0.$$
We will be interested in the cohomology group $H({\cal A}, \delta)$.
A $\bk$-linear functional $\int: {\cal A} \rightarrow \bk$ 
on a DGBV-algebra is called {\em an integral} if
for all $a, b \in {\cal A}$,
\begin{eqnarray}
\int (\delta a) \wedge b & = & 
	(-1)^{|a|+1} \int a \wedge \delta b, \label{int1} \\
\int (\Delta a ) \wedge b & = &
	(-1)^{|a|} \int a \wedge \Delta b. \label{int2}
\end{eqnarray} 
Under these conditions, it is clear that $\int $ induces a
scalar product on $H = H({\cal A}, \delta)$:
$(a, b ) = \int a \wedge b$.
If it is nondegenerate on $H$, we say that the integral
is {\em nice}.
It is obvious that 
$$(\alpha \wedge \beta, \gamma) = (\alpha, \beta \wedge \gamma).$$
Hence if ${\cal A}$ has a nice integral, 
$(H, \wedge, (\cdot, \cdot))$ is a (graded) Frobenius algebra.

Under suitable conditions, one can construct Frobenius manifolds from 
DGBV algebras. The following result is due to Barannikov-Kontsevich \cite{Bar-Kon} 
and Manin \cite{Man2}:

\begin{theorem} \label{thm:construction}
Let $({\cal A}, \wedge, \delta, \Delta, [\cdot \bullet \cdot])$
be a DGBV algebra satisfying the following conditions:
\begin{itemize}
\item[(a)] $H = H({\cal A}, \delta)$ is finite dimensional.
\item[(b)] There is a nice integral on ${\cal A}$.
\item[(c)] The inclusions $(Ker \Delta, \delta) \hookrightarrow 
({\cal A}, \delta)$ and $(Ker \delta, \Delta) \hookrightarrow 
({\cal A}, \Delta)$  induce isomorphisms of cohomology.
\end{itemize}
Then there is a canonical construction of a
formal Frobenius manifold  structure with identity on $H$.
\end{theorem}

We now indicate how to obtain the potential function $\Phi$.
It is based on the existence of a solution 
$\Gamma = \sum \Gamma_n$ to
\begin{eqnarray*}
&\delta \Gamma + \frac{1}{2}[\Gamma \bullet \Gamma] = 0, \\
&\Delta \Gamma = 0, 
\end{eqnarray*}
which satisfies the following conditions:
(a) $\Gamma_0 = 0$;
(b) $\Gamma_1 = \sum x^j e_j$, $e_j \in \Ker \delta \cap \Ker \Delta$,
where the classes of $e_j$'s generates $H = H({\cal A}, \delta)$;
(c) for $n > 1$, $\Gamma_n \in \Img \Delta$ is 
a homogeneous super polynomial of degree $n$ in $x^j$'s,
such that the total degree of $\Gamma_n$ is even; 
(d) $x^0$ only appears in $\Gamma_1$.
Such a solution is called a {\em normalized universal solution}.
Under suitable conditions,
its existence can be established inductively. 
Let $\Gamma = \Gamma_1 + \Delta B$ be a normalized solution,
then
$$\Phi = \int \frac{1}{6}\Gamma^3 - \frac{1}{2}\delta B \Delta B
 = \int \frac{1}{6}\Gamma^3 - \frac{1}{4} \Gamma \wedge \Gamma \wedge
(\Gamma - \Gamma_1).$$

\section{Cartan model of equivariant cohomology}

We will use the Cartan model for equivariant cohomology.  
We refer the readers to Atiyah-Bott \cite{Ati-Bot} and 
Berline-Getzler-Vergne \cite{Ber-Get-Ver} for more details. 
Throughout this paper, $K$ will be a compact connected Lie 
group, with $\fk$ as its Lie algebra.
Let $M$ be a compact smooth $K$-manifold. 
The $K$-action on $M$ induces a homomorphism
from the Lie algebra $\fk$ to the Lie algebra of 
vector fields on $M$.
Let $\{ \xi_a \}$ be a basis of $\fk$, such that
$$[\xi_a, \xi_b] = f_{ab}^c \xi_c, $$
where $f^c_{ab}$'s are the structure constants. 
Let 
$\{ \Theta^a \}$ be the dual basis in $\fk^*$. 
Denote by $\iota_a$ and ${\cal L}_a$ the contraction and 
the Lie derivative by the vector field corresponding 
to $\xi_a \in \fk$ respectively. The Cartan model is given by
the complex $(\Omega_K(X), D_K)$, 
where $\Omega_K(M) =(S(\fk^*) \otimes \Omega(M))^K$, and
$D_K = 1 \otimes d - \Theta^a \otimes \iota_a$, which is called the Cartan differential. 
Since $D_K$ is a $K$-invariant operator on 
$S(\fk^*) \otimes \Omega(M)$, 
it then maps $\Omega_K(M)$ to itself. 
Furthermore,
since $\Theta^a L_a$ acts as zero on $S(\fk^*)$,
we have
$$D_K^2 = - \Theta^a \otimes {\cal L}_a 
= -\Theta^a  (L_a \otimes 1 + 1 \otimes {\cal L}_a).$$ 
Therefore, $D_K^2 = 0$ on $\Omega_K(M)$.
The Cartan model defines 
equivariant cohomology of the $K$-manifold $M$ as
$$H^*_K(M) = \Ker D_K / \Img D_K.$$
The wedge product $\wedge$ on $\Omega^*(M)$ can be 
extended to $\Omega^*_K(M)$.
This makes $\Omega^*_K(M)$ an algebra over $S(\fk^*)^K$.
It is easy to see that $D_K$ is a derivation,
i.e.,
$$D_K(\alpha \wedge \beta) = (D_K\alpha) \wedge \beta
+ (-1)^{|\alpha|}\alpha \wedge D_K\beta,$$
for homogeneous $\alpha, \beta \in \Omega^*_K(M)$.
Hence $H^*_K(M)$ is an algebra over $S(\fk^*)^K$.

Notice that there is a $S(\fk^*)^K$-linear operator 
$$\int_M: \Omega_K^*(M) \to S(\fk^*)^K$$
which is defined by sending  differential forms of degree $\dim (M)$ 
to its integral over $M$, and all other forms to zero.
Since we assume $M$ has no boundary,
by Stokes theorem, 
it is easy to see that
\begin{eqnarray} \label{eqn:Stokes1}
\int_M (D_K\alpha) \wedge \beta 
	= (-1)^{|\alpha|+1}\int_M  \alpha \wedge D_K\beta.
\end{eqnarray}

For simplicity of notation, we will simply write $D_K = d - C$, where 
$C=\Theta^a\iota_a$.
Then $dC + C d =0$, $C^2 = 0$ on $\Omega^*_K(M)$.
There is a natural bigrading on $\Omega^*_K(M)$:
$$(\Omega^*_K(M))^{p, q} = (\Omega^{p-q}(M) \otimes S^q(\fk^*))^K.$$
With respect to this bigrading, $d$ has bidegree $(1, 0)$, 
$C$ has bigrading $(0, 1)$.
Every element $\alpha_K \in \Omega^*_K(M)$ can be written as
$$\alpha_K = \sum_{k \geq 0} \alpha^{(2k)},$$
such that $D_K \alpha_K = 0$ if and only if
$d\alpha^{(0)} = 0$ and
$d\alpha^{(2k+2)} =  C\alpha^{(2k)}$, $k \geq 0$.

\section{Invariant Poisson structure 
and DGBV algebra structure on Cartan model} 

\label{sec:Poisson}

We now assume that $M$ has an $K$-invariant Poisson structure $w$, 
i.e. $w \in \Gamma(M, \Lambda^2 TM)^K$ and 
the Schouten-Nijenhuis bracket $[w, w] = 0$.
For any Poisson structure $w$,
Koszul \cite{Kos} defined an operator 
$\Delta: \Omega^*(M) \to \Omega^{*-1}(M)$ by
$\Delta = [\iota_w, d]$.
He also showed that $\Delta$ has the following important properties:
\begin{eqnarray*}
\Delta^2 = 0, & [d, \Delta] = d\Delta + \Delta d = 0,
\end{eqnarray*}
and if we set
$$[\alpha, \beta]_{\Delta} = (-1)^{|\alpha|}(\Delta(\alpha \wedge \beta) 
	- (\Delta \alpha) \wedge \beta
	- (-1)^{|\Delta|} \alpha \wedge \Delta \beta),$$
then 
$$[\alpha, \beta \wedge \gamma]_{\Delta} 
= [\alpha, \beta \wedge \gamma]_{\Delta} \wedge \gamma
+(-1)^{(|\alpha|+1)|\beta|} \beta \wedge [\alpha, \gamma]_{\Delta}.$$
This implies that 
$(\Omega^*(M), \wedge, d, \Delta, [\cdot, \cdot]_{\Delta})$ is a DGBV algebra.
In our case, we can extend $\Delta$ and $[\cdot, \cdot]_{\Delta}$ 
to $\Omega^*(M)\otimes S(\fk^*)$.
It is clear that they both commute with the group action.
Hence they restrict to $\Omega^*_K(M)$.

\begin{proposition}
Let $M$ be a $K$-manifold with a $K$-invariant Poisson structure $w$,
then $(\Omega^*_K(M), \wedge, D_K, \Delta, [\cdot, \cdot]_{\Delta})$
is a DGBV algebra.
\end{proposition}

\begin{proof}
One only needs to prove $[D_K, \Delta] = 0$.
Now $D_K = d - C = d - \Theta^a\iota_a$.
We have $[d, \Delta] = 0$,
and
\begin{eqnarray*}
[\iota_a, \Delta] = [\iota_a, [\iota_w, d]]
= [[\iota_a, \iota_w], d] + [\iota_w, [\iota_a,  d]] 
= [\iota_w, \cL_a] = -[\iota_w, L_a]= 0
\end{eqnarray*}
on $\Omega^*_K(M)$, hence $[C, \Delta] = 0$.
The proof is complete.
\end{proof}

From \cite{Cao-Zho2}, we also have 
\begin{lemma}
For a $K$-manifold $M$ with a $K$-invariant Poisson structure $w$,
we have
\begin{eqnarray*}
\int_M (\Delta \alpha) \wedge \beta = 
	(-1)^{|\alpha|} \int \alpha \wedge \Delta \beta.
\end{eqnarray*}
\end{lemma}

\section{Symplectic manifolds with Hamiltonian actions}

We now assume that $M$ has a symplectic structure $\omega$
and the $K$-action is Hamiltonian,
i.e., the $K$-action preserves $\omega$ and
there is a $K$-equivariant map $\mu: M \to \fk^*$,
such that
$$d\langle \mu, \xi_a \rangle = \iota_a \omega.$$
The symplectic structure $\omega$ induces an isomorphism $T^*M \cong TM$,
hence isomorphisms $\Omega^*(M) \cong \Gamma(M, \Lambda^*TM)$.
Denote by $w$ the bi-vector field corresponding to $\omega$.
Then $w$ is an invariant Poisson structure. Hence, we have

\begin{proposition}
For a symplectic manifold $M$ with a Hamiltonian $K$-action,
$(\Omega^*_K(M), \wedge, D_K, \Delta, [\cdot, \cdot]_{\Delta})$
is a DGBV algebra.
\end{proposition}

\begin{remark}
It is tempting to define $\Delta_K = \Delta - d\mu \wedge$. Indeed,
it is easy to show that $\Delta_K^2 = [D_K, \Delta_K] = 0$.
However, we do not have 
$$[\alpha, \beta \wedge \gamma]_{\Delta_K} 
= [\alpha, \beta \wedge \gamma]_{\Delta_K} \wedge \gamma
+(-1)^{(|\alpha|+1)|\beta|} \beta \wedge [\alpha, \gamma]_{\Delta_K}.$$
Hence  
$(\Omega^*_K(M), \wedge, D_K, \Delta_K, [\cdot, \cdot]_{\Delta_K})$
is not a DGBV algebra.
\end{remark}



For a closed symplectic manifold $M$ with a Hamiltonian $K$-action,
a result of Kirwan \cite{Kir} (p. 68, Proposition 5.8) says that 
$H^*_K(M) \cong H^*(M) \otimes_{\bR} S(\fk^*)^K$ as vector spaces over $\bR$.
An important consequence of the above result of Kirwan is that
every de Rham cohomology class of $M$ has a representative $\alpha$,
which can be extended to a $D_K$ closed form $\alpha_K$ of the form
$$\alpha_K=\alpha + \Theta^a \alpha_a + \cdots.$$
Therefore,
one can find $D_K$-closed forms 
$\{\alpha_{Ki} =\alpha_i +\Theta^a \alpha_{ia} + \cdots \}$
such that the matrix $(\int_M \alpha_i \wedge \alpha_j)$ is invertible over
$\bR$. However the matrix $(\int_M \alpha_{Ki} \wedge \alpha_{Kj})$ may not be invertible
over $S(\fk^*)^K$. 
Later we will prove that for Hamiltonian actions on a closed K\"{a}hler 
manifold by holomorphic isometries,
we can find natural extensions $\{\alpha_{Ki}\}$ such that 
$$\int_M \alpha_{Ki} \wedge \alpha_{Kj} = \int_M \alpha_i \wedge \alpha_j.$$
For now, to invert the matrix  $(\int_M \alpha_{Ki} \wedge \alpha_{Kj})$, 
we need to work over a field.
Denote by $T$ a maximal torus of $K$,
$\ft$ its Lie algebra and $W$ the Weyl group.
Then $S(\fk^*)^K = S(\ft^*)^W$.
Hence $S(\fk^*)^K$ is an integral domain, since $S(\ft^*)$ is a polynomial algebra.
Denote by $F(\fk^*)$ its fractional field,
i.e. $F(\fk^*) = \{f/g: f, g \in S(\fk^*)^K\}$.
Define
$$\widetilde{\Omega}^*_K(M) = \Omega^*_K(M) \otimes_{S(\fk^*)^K} F(\fk^*),$$
Extend $D_K$, $\wedge$, $\Delta$ etc. to $\widetilde{\Omega}^*_K(M)$
and define
$$\widetilde{H}^*_K(M) = H^*(\widetilde{\Omega}^*_K(M), D_K).$$
Then we have 
$$\widetilde{H}^*_K(M) = H^*_K(M) \otimes_{S(\fk^*)^K} F(\fk^*)$$
as vector spaces over $F(\fk^*)$.
Now the matrix $(\int_M \alpha_{Ki} \wedge \alpha_{Kj})$ 
has a nonzero determinant,
hence it is invertible over $F(\fk^*)$. 
Therefore  
$(\widetilde{\Omega}^*_K(M), \wedge, D_K, \Delta, [\cdot, \cdot]_{\Delta})$
satisfies conditions (a) and (b) in Theorem \ref{thm:construction}
over the field $\bk = F(\fk^*)$.
Thus we have the following 

\begin{theorem}
Let $M$ be a closed symplectic manifold with a Hamiltonian $K$-action.
Suppose that the inclusions $i: (\Ker \Delta, D_K) \hookrightarrow
(\widetilde{\Omega}^*_K(M), D_K)$
and $j: (\Ker D_K, \Delta) \hookrightarrow (\widetilde{\Omega}^*_K(M), \Delta)$
induce isomorphisms on cohomology.
Then over the field $F(\fk^*)$, 
the DGBV algebra 
$(\widetilde{\Omega}^*_K(M), \wedge, D_K, \Delta, [\cdot, \cdot]_{\Delta})$
satisfies all the conditions in 
Theorem \ref{thm:construction}.
Hence there is a canonical construction of formal Frobenius manifold
structure on $\widetilde{H}^*_K(M)$.
\end{theorem}

\section{K\"{a}hler manifolds with holomorphic Hamiltonian actions}

We now further restrict our attention to a closed K\"{a}hler manifold $M$,
such that $K$ acts on $M$ by holomorphic isometries.
Then the K\"{a}hler form is an invariant symplectic form,
hence the results in \S \ref{sec:Poisson} apply.
The main advantage here is that for a K\"{a}hler manifold, 
we can exploit some nice features of the Hodge theory to establish the quasi-isomorphisms
property (c) in Theorem \ref{thm:construction}. Note that in Lemma 5.2-5.4, we will not require
the $K$-action to be Hamiltonian.

The almost complex structure $J: TM \to TM$ induces a decomposition
$TM \otimes_{\bR} \bC = T^{1, 0}M \oplus T^{0, 1}M$.
There is an induced decomposition
$\Omega^*(M) \otimes_{\bR} \bC = \Omega^{*, *}(M)$,
and $d = \partial + \bar{\partial}$,
where
\begin{eqnarray*}
\partial: \Omega^{*, *}(M) \to \Omega^{*+1, *}(M), & 
\bar{\partial}: \Omega^{*, *}(M) \to \Omega^{*, *+1}(M).
\end{eqnarray*}
Denote also by $J$ the linear map on $\Omega^{*, *}(M)$ induced by $J$.
Then we have
$$J\alpha = (-1)^q i^{p+q} \alpha,$$
for $\alpha \in \Omega^{p, q}(M)$.
Hence $J^2 = (-1)^{p+q}$ and $J^{-1} = (-1)^p i^{p+q}$ 
on $\Omega^{p, q}(M)$.
It is easy to see that
\begin{eqnarray*}
J^{-1} \partial J = - i \partial, & 
J^{-1} \bar{\partial} J = i \bar{\partial}.
\end{eqnarray*}
Hence
$$J^{-1} d J = J^{-1}\partial J + J^{-1} \bar{\partial} J 
= i( \bar{\partial} - \partial) = d^c.$$
We are also interested in the formal adjoints 
$d^*$, $\partial^*$ and $\bar{\partial}^*$.
Since
\begin{eqnarray*}
\partial^*: \Omega^{*, *}(M) \to \Omega^{*-1, *}(M), & 
\bar{\partial}^*: \Omega^{*, *}(M) \to \Omega^{*, *-1}(M).
\end{eqnarray*}
Hence
\begin{eqnarray*}
J^{-1} \partial^* J =  i \partial^*, & 
J^{-1} \bar{\partial}^* J = -i \bar{\partial}^*.
\end{eqnarray*}
Therefore
$$J^{-1} d^* J = J^{-1}\partial^* J + J^{-1} \bar{\partial}^* J 
= -i(\bar{\partial}^* - \partial^*) = (d^c)^*.$$
With the help of K\"{a}hler identities,
one can obtain the following well-known equalities
(see e.g. Deligne-Griffiths-Morgan-Sullivan \cite{Del-Gri-Mor-Sul}):
\begin{eqnarray*}
&& d^2 = (d^c)^2 = dd^c + d^cd = 0, \\
&& (d^*)^2 = ((d^c)^*)^2 = d^*(d^c)^* + (d^c)^*d^* = 0, \\
&& d(d^c)^* + (d^c)^*d = d^*d^c + d^cd^* = 0, \\
&& dd^* + d^*d = d^c(d^c)^* + (d^c)^* d^c = \square.
\end{eqnarray*}
Here, $\square$ in the last equality denotes the Laplace operator on forms.  
As a consequence, one has the following
Hodge decompositions 
\begin{eqnarray*}
\Omega^*(M) & = & \cH \oplus \Img d \oplus \Img d^* = 
\cH \oplus \Img d^c \oplus \Img (d^c)^* \\
& = & \cH \oplus \Img dd^c \oplus \Img d^*d^c \oplus \Img d(d^c)^* 
\oplus \Img d^*(d^c)^*,
\end{eqnarray*}
where $\cH$ is the space of harmonic forms.

\begin{lemma} \label{lm:Key}
On a closed K\"{a}hler manifold $M$,
if $\Delta d\beta = 0$ for some $\beta \in \Omega_K^*(M)$,
then there exist $\beta^H \in \Ker \square$, $a, b, c \in \Omega_K^*(M)$,
such that
$$\beta = \beta^H + \Delta^*d a + \Delta d b + \Delta d^* c.$$
\end{lemma}

\begin{proof}
It suffices to prove the result for $\Omega^*(M)$.
The extension to $\Omega^*_K(M)$ is straightforward. But for $\Omega^*(M)$,
K\"{a}kler identity implies that $\Delta = -(d^c)^*$,
so the lemma follows from the above five-fold decomposition.
\end{proof}

\begin{lemma} \label{lm:jinj}
The inclusion 
$j: (\Ker D_K, \Delta) \hookrightarrow (\Omega^*_K(M), \Delta)$ 
induces an injective map in cohomology.
\end{lemma}

\begin{proof}
We need to show that if $D_K \alpha_K = 0$ and $\alpha_K = \Delta \beta_K$
for some $\beta_K \in \Omega^*_K(M)$,
then there exists $\beta_K' \in \Ker D_K$, 
such that $\alpha_K = \Delta \beta'_K$.
We use the bigrading on $\Omega^*_K(M)$ to write
\begin{eqnarray*}
\alpha_K = \sum_{k\geq 0} \alpha^{(2k)}, &
\beta_K =  \sum_{k\geq 0} \beta^{(2k)},
\end{eqnarray*}
such that 
\begin{align*} 
d\alpha^{(0)} &= 0, & \alpha^{(0)} &= \Delta \beta^{(0)}, \\
d\alpha^{(2)} & = C \alpha^{(0)}, & \alpha^{(2)} &= \Delta \beta^{(2)}, \\
\cdots\cdots &
\end{align*}
We will repeatedly use the following corollary of Lemma \ref{lm:Key}: 
if $d\Delta \beta = 0$ for some $\beta \in \Omega^*_K(M)$,
then there exists $\gamma \in \Omega^*_K(M)$,
such that
$$\Delta \beta = \Delta d\gamma.$$
In fact, we can take $\gamma = Gd^*\beta$, where $G$ is the Green operator, so that
\begin{eqnarray*}
\Delta d\gamma =\Delta dGd^*\beta = Gdd^*\Delta \beta = G(\square - d^*d) \Delta \beta 
= G\square \Delta \beta = \Delta \beta.
\end{eqnarray*}
Now $d\Delta \beta^{(0)} = d\alpha^{(0)} = 0$,
hence 
$$\alpha^{(0)} = \Delta \beta^{(0)} = \Delta d \gamma^{(0)},$$
where $\gamma^{(0)} = Gd^*\beta^{(0)}$.
Also, 
\begin{eqnarray*}
d\Delta\beta^{(2)} = d\alpha^{(2)} = C\alpha^{(0)} = C\Delta d \gamma^{(0)}
= -d\Delta C\gamma^{(0)}.
\end{eqnarray*}
Hence $\Delta d(\beta^{(2)} + C\gamma^{(0)}) = 0$.
Therefore,
$$\Delta(\beta^{(2)} + C\gamma^{(0)}) = \Delta d\gamma^{(2)},$$ 
where $\gamma^{(2)} = Gd^*(\beta^{(2)} + C\gamma^{(0)})$.
Equivalently, we have
$$\alpha^{(2)} = \Delta\beta^{(2)} 
=  \Delta d \gamma^{(2)} - \Delta C\gamma^{(0)}
= \Delta (d\gamma^{(2)} - C\gamma^{(0)}).$$
Inductively, we have for $k \geq 0$,
$$\alpha^{(2k+2)} = \Delta (d \gamma^{(2k+2)} - C \gamma^{(2k)}),$$
where $\gamma^{(2k+2)} = Gd^* (\beta^{(2k+2)} + C \gamma^{(2k)})$.
Setting $\beta'_K= D_K \sum_{k\geq 0} \gamma^{(2k)}$, it is then clear that
$D_K\beta'_K=0$ and 
$$\alpha_K = \Delta \beta_K = \Delta \beta'_K.$$
\end{proof}

\begin{remark}
We actually prove the following stronger result:
$$\Ker D_K \cap \Img \Delta = \Img \Delta D_K.$$
\end{remark}
 
\begin{lemma} \label{lm:iinj}
The inclusion $i: (\Ker \Delta, D_K) \hookrightarrow (\Omega^*_K(M), D_K)$
induces an injective map in cohomology.
\end{lemma}

\begin{proof}
We need to show that if $\Delta \alpha_K = 0$ and $\alpha_K = D_K \beta_K$
for some $\beta_K \in \Omega^*_K(M)$,
then there exists $\beta'_K \in \Ker \Delta$,
such that  $\alpha_K = D_K \beta_K'$.
We will use repeatedly the following corollary of Lemma \ref{lm:Key}:
if $d\Delta \beta = 0$ for some $\beta \in \Omega_K^*(M)$,
then there exist $\beta^H \in \Ker \square$, $\phi, \psi \in \Omega^*_K(M)$,
such that
$$\beta = \beta^H + \Delta \phi + d\psi.$$
Decompose $\alpha_K$ and $\beta_K$ as usual.
We have
\begin{align*}
\Delta \alpha^{(0)} & = 0, & \alpha^{(0)} & = d \beta^{(0)}, \\
\Delta \alpha^{(2)} & = 0, 
	& \alpha^{(2)} & = d \beta^{(2)} - C \beta^{(0)}, \\
\cdots\cdots&
\end{align*}
Now $\Delta d\beta^{(0)} = \Delta \alpha^{(0)} = 0$,
hence
$$\beta^{(0)} = (\beta^{(0)})^H + \Delta \phi^{(0)} + d \psi^{(0)}.$$
Therefore,
\begin{eqnarray*}
 \Delta d \beta^{(2)} = \Delta (\alpha^{(2)} + C \beta^{(0)})
= - C\Delta \beta^{(0)}
= -C\Delta d \psi^{(0)} = -\Delta d C\psi^{(0)}.
\end{eqnarray*}
So we have $\Delta d(\beta^{(2)} + C\psi^{(0)}) = 0$,
hence
$$\beta^{(2)} + C\psi^{(0)} 
= (\beta^{(2)} + C\psi^{(0)})^H + \Delta \phi^{(2)} + d \psi^{(2)}.$$
By induction, we find that for $k \geq 0$,
$$\beta^{(2k+2)} = (\beta^{(2k+2)} + C\psi^{(2k)})^H 
+ \Delta \phi^{(2k+2)} + d \psi^{(2k+2)} - C \psi^{(2k)}.$$
Setting $\phi_K = \sum_{k \geq 0} \phi^{(2k)}$ 
and $\psi_K = \sum_{k \geq 0} \psi^{(2k)}$,
then we have
$$\beta_K = (\beta_K + C \psi_K)^H + \Delta \phi_K + D_K\psi_K.$$
Hence
\begin{eqnarray} \label{eqn:??}
\alpha_K = D_K \beta_K = D_K((\beta_K + C \psi_K)^H + \Delta \phi_K).
\end{eqnarray}
Since $\Delta ((\beta_K + C \psi_K)^H + \Delta \phi_K)=0$,
the proof is complete.
\end{proof}

\begin{remark}
We actually have proved $\Ker \Delta \cap \Img D_K = D_K \Ker \Delta$.
If one can show that for any $\beta \in \Ker \square$,
there exists $\gamma \in \Omega^*_K(M)$,
such that $D_K (\beta + \Delta \gamma) = 0$,
then by (\ref{eqn:??}),
$\Ker \Delta \cap \Img D_K = \Img D_K\Delta$.
\end{remark}

\begin{lemma}  \label{lm:isur}
The inclusion $i: (\Ker \Delta, D_K) \hookrightarrow (\Omega^*_K(M), D_K)$
induces a surjective map in cohomology.
\end{lemma}

\begin{proof}
We need to show that 
If $D_K \alpha_K = 0$ for some $\alpha_K\in \Omega^*_K(M)$,
then there exists $\beta_K \in \Omega^*_K(M)$,
such that $\Delta(\alpha_K - D_K\beta_K) = 0$.
Decompose $\alpha_K$ as usual,
then we have
$d\alpha^{(0)} = 0$, $d\alpha^{(2)} = C\alpha^{(0)}$, etc.
First of all,
there exists $\beta^{(0)} \in \Omega^*_K(M)$,
such that $\alpha^{(0)} - d \beta^{(0)} = (\alpha^{(0)})^H \in \Ker \square$.
Now
\begin{eqnarray*}
 \Delta d(\alpha^{(2)} + C \beta^{(0)}) 
= \Delta C \alpha^{(0)} - \Delta C d\beta^{(0)} =\Delta C (\alpha^{(0)})^H = 0.
\end{eqnarray*}
Hence $\alpha^{(2)} + C \beta^{(0)} 
= (\alpha^{(2)} + C \beta^{(0)})^H + \Delta \gamma^{(2)} + d\beta^{(2)}$.
By induction,
we can find $\beta^{(2k+2)} \in \Omega^*_K(M)$, for $k \geq 0$,
such that
\begin{eqnarray*}
\alpha^{(2k+2)} 
= (\alpha^{(2k+2)} + C \beta^{(2k)})^H + \Delta \gamma^{(2k+2)} 
+ d\beta^{(2k+2)} - C \beta^{(2k)}.
\end{eqnarray*}
Set $\beta_K = \sum_{k \geq 0} \beta^{(2k)}$, 
$\gamma_K = \sum_{k \geq 1} \gamma^{(2k)}$.
Then we have
$$\alpha_K = (\alpha_K + C \beta_K)^H + \Delta \gamma_K + D_K\beta_K.$$
Hence 
$$\Delta(\alpha_K - D_K \beta_K) 
= \Delta((\alpha_K + C \beta_K)^H + \Delta \gamma_K) = 0.$$
\end{proof}

From now on we assume that the $K$-action is Hamiltonian,
and let $\mu = \Theta^a\mu_a$ be the moment map.

\begin{lemma} \label{lm:contraction}
$\iota_a \alpha = \mu_a \Delta \alpha - \Delta (\mu_a \alpha)$.
\end{lemma}

\begin{proof}
For $\alpha, \beta \in \Omega^*(M)$, we have
\begin{eqnarray*}
&& \langle \iota_a \alpha, \beta \rangle
= \langle \alpha, Jd\mu_a \wedge \beta \rangle  
= \langle J \alpha, J(Jd\mu_a \wedge \beta) \rangle \\
& = & - \langle  J \alpha, d\mu_a \wedge J\beta \rangle  
= - \langle J \alpha, d (\mu_a J\beta) - \mu_a dJ\beta \rangle \\
& = & - \langle d^* J \alpha, \mu_a J\beta \rangle   
	+ \langle \mu_aJ \alpha, dJ\beta\rangle \\
& = & - \langle \mu_a d^* J \alpha,  J\beta \rangle   
	+ \langle d^*J(\mu_a \alpha), J\beta \rangle \\
& = & - \langle J^{-1}(\mu_ad^* J \alpha) - J^{-1}d^*J(\mu_a \alpha), 
	\beta\rangle   \\
& = & - \langle \mu_a(d^c)^*\alpha) - (d^c)^*(\mu_a \alpha), \beta\rangle \\
& = & \langle \mu_a \Delta \alpha - \Delta (\mu_a \alpha), \beta \rangle.
\end{eqnarray*}
\end{proof}

\begin{corollary} \label{cor:contraction}
If $\Delta\alpha = 0$,
then $C(\alpha) = -\Delta (\mu \alpha)$.
\end{corollary}

\begin{proposition} \label{prop:harmonic}
On a closed K\"{a}hler manifold $M$ with a Hamiltonian $K$-action by
holomorphic isometries,
any harmonic form 
$\alpha^{(0)}$ can be canonically extended to a $D_K$-closed form
$\sum_{k \geq 0} \alpha^{(2k)}$,
where $\alpha^{(2k+2)} = - Gd^*\Delta (\mu\alpha^{(2k)}) 
\in \Img \Delta d^*$. 
\end{proposition}

\begin{proof}
First notice that $\alpha^{(0)} \in \Omega^*(M)^K$.
Indeed, since $K$ is connected,
the action of $K$ on $H^*(M)$ is trivial.
Hence for any $g \in G$,
$g(\alpha^{(0)})$ is a harmonic form 
in the same cohomology class as $\alpha^{(0)}$,
therefore, $g(\alpha^{(0)}) = \alpha^{(0)}$.
By Corollary \ref{cor:contraction},
$C(\alpha^{(0)}) = -\Delta (\mu\alpha^{(0)})$.
Since $d  C(\alpha^{(0)}) = -C d\alpha^{(0)} = 0$ and 
 $C(\alpha^{(0)})$ has no harmonic part, we have 
$$C(\alpha^{(0)})=d \alpha^{(2)}$$
where $\alpha^{(2)} =-Gd^*\Delta (\mu\alpha^{(0)})=Gd^* C(\alpha^{(0)})\in \Omega^*_K(M)$.
Now since $\alpha^{(2)} \in \Img \Delta$, 
$C(\alpha^{(2)}) = - \Delta(\mu\alpha^{(2)})$.
From 
$$dC(\alpha^{(2)}) = -Cd\alpha^{(2)} = - C^2(\alpha^{(0)}) = 0,$$
we see that if we set 
$\alpha^{(4)} = -Gd^*\Delta (\mu \alpha^{(2)})$,
then $\alpha^{(4)} \in \Omega^*_K(M)$ and
we have $d \alpha^{(4)} = C(\alpha^{(2)})$.
Inductively, one gets $\alpha^{(2k)}$
in the same way. 
This process terminates after finitely many steps,
since each time the degree of the differential forms are reduced by $2$.
Then $\sum_{k \geq 0} \alpha^{(2k)}$ is a $D_K$-closed form.
\end{proof}

As a corollary, 
we get an easy proof of  
$$H^*_K(M) \cong H^*(M) \otimes S(\fk^*)^K$$
in the case of closed K\"{a}hler manifolds.
As another corollary, we get the following

\begin{lemma} \label{lm:jsur}
On a closed K\"{a}hler manifold $M$ with a Hamiltonian $K$-action by
holomorphic isometries,
the inclusion $j: (\Ker D_K, \Delta) 
\hookrightarrow (\Omega^*_K(M), \Delta)$
induces a surjective map in cohomology.
\end{lemma}

\begin{proof}
It suffices to show that for any $\alpha^{(0)} \in \Ker \Delta$,
there exists $\beta \in \Omega^*_K(M)$,
such that $D_K(\alpha^{(0)} + \Delta \beta) = 0$.
Without loss of generality,
we can assume that $\alpha^{(0)}$ is harmonic.
By Proposition \ref{prop:harmonic}, 
we can take 
$$\beta = Gd^*\sum_{k \geq 0} \mu \alpha^{(2k)},$$
where $\alpha^{(2k+2)} = - Gd^*\Delta (\mu\alpha^{(2k)})$.
\end{proof}

Combining Lemmas \ref{lm:jinj}, \ref{lm:iinj}, \ref{lm:isur} and 
\ref{lm:jsur},
we get

\begin{theorem}
On a closed K\"{a}hler manifold $M$ with a Hamiltonian $K$-action by
holomorphic isometries,
the inclusions
$i: (\Ker \Delta, D_K) \hookrightarrow (\Omega^*_K(M), D_K)$ and
$j: (\Ker D_K, \Delta) \hookrightarrow (\Omega^*_K(M), \Delta)$
induce isomorphisms in cohomology.
\end{theorem}

\section{Normalized universal solution and formal Frobenius manifold structure}

Assume that $\{ \omega_a^{(0)} \in \cH\}$ gives rise to 
a homogeneous basis of $H^*(M)$.
By Proposition \ref{prop:harmonic},
each $\omega^{(0)}_a$ can be extended to a $D_K$-closed class
$$\omega_{Ka} = \sum_{k \geq 0} \omega_a^{(2k)},$$
where  $\alpha_a^{(2k+2)} = - Gd^*\Delta (\mu\alpha_a^{(2k)})$ for $k \geq 0$.
In particular,
$\alpha_{K0} = 1$.
Now $\{ \omega_{Ka} \}$ are free generators of $H^*_K(M)$.
Write $\omega_{Ka} = \omega_a^{(0)} + \Delta \gamma_a$,
then from Lemma 3.1 we have
\begin{eqnarray*}
&& \int_M \omega_{Ka} \wedge \omega_{Kb} 
= \int_M (\omega_a^{(0)} + \Delta \gamma_a) \wedge 
	(\omega_b^{(0)} + \Delta \gamma_b) \\
& = & \int_M \omega_a^{(0)}\wedge \omega_b^{(0)} + 
	\Delta \gamma_a \wedge \omega_b^{(0)} 
	+ \omega_a^{(0)}  \wedge \Delta \gamma_b
	+ \Delta \gamma_a \wedge\Delta \gamma_b \\
& = & \int_M \omega_a^{(0)}\wedge \omega_b^{(0)} \pm 
	\gamma_a \wedge \Delta \omega_b^{(0)} 
	\pm \Delta \omega_a^{(0)}  \wedge \gamma_b
	\pm\Delta^2 \gamma_a \wedge \gamma_b \\
& = & \int_M \omega_a^{(0)}\wedge \omega_b^{(0)}.
\end{eqnarray*}
In other words,
the matrix $(\eta_{ab})$ for the pairing $(\cdot, \cdot)$ 
is the same as in the ordinary case.
Hence we can take $\bk = S(\fk^*)^K$ in Theorem \ref{thm:construction}. Thus, we have

\begin{theorem}
Let $M$ be a closed K\"{a}hler manifold with a Hamiltonian $K$-action by
holomorphic isometries,
then there is a formal Frobenius manifold structure 
on $H^*_K(M)$ obtained by DGBV algebraic construction over $S(\fk^*)^K$.
\end{theorem}

To obtain the normalized universal solution, we take $\Gamma_{K1} = x^a\alpha_{Ka}$.
For $n > 1$, we find $\Gamma_{Kn} \in \Img \Delta$
by inductively solving 
$$D_K\Gamma_{Kn} 
= -\frac{1}{2}\sum_{p=1}^{n-1}[\Gamma_{Kp} \bullet \Gamma_{Kn-p}]_{\Delta}
= - \frac{1}{2} \sum_{p=1}^{n-1}\Delta(\Gamma_{Kp} \wedge \Gamma_{Kn-p}).$$
Standard argument shows that the right hand side is $D_K$-closed,
hence by the proof of Lemma \ref{lm:jinj},
we take
$$\Gamma_{Kn} 
= \frac{1}{2} \Delta \sum_{k \geq 0} \gamma_n^{(2k)},$$
where we set
\begin{eqnarray*}
&&\beta_n = \sum_{p =1}^{n-1}(\Gamma_p \wedge \Gamma_{n-p}), \\
&&\gamma_n^{(0)} = Gd^*\beta^{(0)}_n, \\
&&\gamma_n^{(2k+2)} = Gd^*(\beta_n^{(2k+2)} + C\gamma_n^{(2k)}), 
	\quad k \geq 0.
\end{eqnarray*}
We have $\Delta \gamma_n^{(0)} = G\Delta d^*\beta^{(0)}_n$.
By Lemma \ref{lm:contraction},
$C\gamma^{(2k)} = \mu \Delta \gamma^{(2k)} - \Delta (\mu \gamma^{(2k)})$.
Hence for $k \geq 0$,
$$\Delta \gamma_n^{(2k+2)} = 
G\Delta d^*(\beta_n^{(2k+2)}  + \mu \Delta \gamma_n^{(2k)}).$$
Set $\phi_n^{(2k)} = \Delta \gamma_n^{(2k)}$,
then 
$$\Gamma_{Kn} = \frac{1}{2}\sum_{k \geq 0} \phi_n^{(2k)},$$
where
\begin{eqnarray*}
&& \phi_n^{(0)} = G\Delta d^* \beta_n^{(0)}, \\
&& \phi_n^{(2k+2)} = G\Delta d^*(\beta_n^{(2k+2)} + \mu\phi_n^{(2k)}),
	\quad k \geq 0.
\end{eqnarray*}
The potential function
$$\Phi_K = \int_M \frac{1}{6}\Gamma_K^3 - \frac{1}{4} \Gamma_K^2 \wedge
(\Gamma_K - \Gamma_{K1})$$
is a formal power series with coefficients in  $S(\fk^*)^K$.
There is a ring homomorphism $f: S(\fk^*)^K \to S^0(\fk^*)^K = \bR$.
Now $f(\Gamma_K)$ and $f(\Phi_K)$ give exactly the formal Frobenius manifold
structure constructed in Cao-Zhou \cite{Cao-Zho2}.
Hence $\Phi_K$ should be thought of as a family 
of formal Frobenius manifold structures.
This is more clearly seen when $K$ is a torus $T$.
There is a deformation family of associative algebraic structures on $H^*(M)$
given by
$$\int_M \frac{1}{6}\Gamma_{T1}^3.$$
$\Phi_T$ then gives a family of solutions to WDVV equations
with them as initial data.

\end{document}